\documentclass{amsart}

\usepackage{fancyhdr}
\usepackage{lastpage}
\usepackage{stmaryrd,yhmath}

\newcommand{\bdis}{\begin{displaymath}}
\newcommand{\edis}{\end{displaymath}}
\newcommand{\be}{\begin{equation}}
\newcommand{\ee}{\end{equation}}
\newcommand{\mbb}{\mathbb}
\newcommand{\mcal}{\mathcal}

\newcommand{\vp}{\varphi}

\newcommand{\vth}{\vartheta}

\newcommand{\zf}{\zeta\left(\frac{1}{2}+it\right)}

\theoremstyle{definition}

\theoremstyle{remark}
\newtheorem{remark}[]{Remark}

\newtheorem*{mydef1}{{\bf Theorem}}

\newtheorem*{mydef7}{{\bf Question}}

\newtheorem*{mydefSF}{{\bf Spectral Formula}}

\numberwithin{equation}{section}

\begin{document}

\title{Jacob's ladders, Riemann's oscillators, quotient of two oscillating multiforms and set of metamorphoses of
this system}

\author{Jan Moser}

\address{Department of Mathematical Analysis and Numerical Mathematics, Comenius University, Mlynska Dolina M105, 842 48 Bratislava, SLOVAKIA}

\email{jan.mozer@fmph.uniba.sk}


\begin{abstract}
In this paper we introduce complicated oscillating system, namely quotient of two multiforms based on Riemann-Siegel formula. We prove
that there is an infinite set of metamorphoses of this system (=chrysalis) on critical line $\sigma=\frac 12$ into a butterfly
(=infinite series of M\" obius functions in the region of absolute convergence $\sigma>1$).

To memory of the Hardy's \emph{Pure Mathematics}
\end{abstract}
\maketitle

\section{Introduction}

\subsection{}

Let us remind the Riemann-Siegel formula
\be \label{1.1} \begin{split}
& Z(t)=2\sum_{n\leq\tau(t)}\frac{1}{\sqrt{n}}\cos\{\vth(t)-t\ln n\}+\mcal{O}(t^{-1/4}), \\
& Z(t)=e^{i\vth(t)}\zf,\quad \tau(t)=\sqrt{\frac{t}{2\pi}},
\end{split}
\ee
(see \cite{7}, p. 60, comp. \cite{8}, p. 79), where
\bdis
\vth(t)=-\frac t2\ln\pi+\mbox{Im}\-\ln\Gamma\left(\frac 14+i\frac t2\right),
\edis
(see \cite{8}, p. 239). Next, we have defined in our paper \cite{6} the following multiform
\be \label{1.2}
\begin{split}
& G(x_1,\dots,x_k)=\prod_{r=1}^k |Z(x_r)|=\prod_{r=1}^k\left|\zeta\left(\frac 12+ix_r\right)\right|, \\
& x_r>T>0,\ k=1,2,\dots,k_0,\quad k_0\in \mbb{N}
\end{split}
\ee
in connection with the Riemann-Siegel formula (\ref{1.1}). Furthermore, we have defined the subset of the set of points
\bdis
\begin{split}
& (x_1,\dots,x_k)\in [T,+\infty)^k:\ T<x_1<x_2<\dots<x_k, \\
& x_r\not=\gamma:\ \zeta\left(\frac 12+i\gamma\right)=0,\ r=1,\dots,k,
\end{split}
\edis
and we have proposed the following

\begin{mydef7}
Is there in this subset a point
\bdis
(y_1,\dots,y_k)
\edis
of metamorphosis of the multiform (\ref{1.2}) that is also the point of significant change of the structure of the multiform (\ref{1.2})?
\end{mydef7}

\subsection{}

We have obtained the following answer (see \cite{6}): there is an infinite set of elements
\be \label{1.3}
\{ \alpha_0(T),\alpha_1(T),\dots,\alpha_k(T)\},\ T\in (T_0,+\infty),\ T_0>0,
\ee
where $T_0$ is sufficiently big, such that
\be \label{1.4}
\begin{split}
& \prod_{r=1}^k\left|\sum_{n\leq \tau(\alpha_r)}\frac{2}{\sqrt{n}}\cos\{\vth(\alpha_r)-\alpha_r\ln n\}+\mcal{O}(\alpha_r^{-1/4})\right|\sim \\
& \sim \sqrt{\frac{\Lambda}{\left|\sum_{n\leq \tau(\alpha_0)}\frac{2}{\sqrt{n}}\cos\{\vth(\alpha_0)-\alpha_0\ln n\}+\mcal{O}(\alpha_0^{-1/4})\right|}}, \\
& T\to\infty ,
\end{split}
\ee
where (see \cite{6}, (1.7))
\bdis
\Lambda=\sqrt{2\pi}\frac{\sqrt{H}}{H_k}\ln^kT,
\edis
i.e. to the infinite subset
\bdis
\{\alpha_1(T),\dots,\alpha_k(T)\},\ T\in (T_0,+\infty)
\edis
an infinite set of metamorphoses of the multiform (\ref{1.2}) into quite distinct form on the right-hand side of (\ref{1.4}) corresponds.

\begin{remark}
We shall call the elements of the set (\ref{1.3}) as a control parameters (functions) of the metamorphosis. The reason to do so is that the parameters
\bdis
\{\alpha_1(T),\dots,\alpha_k(T)\}
\edis
change the old form into the new one (see (\ref{1.4}) and this last is controlled by the parameter $\alpha_0(T)$. That is, the set (\ref{1.3})
plays a similar role as the shem-ha-m'forash in Golem's metamorphosis.
\end{remark}

\subsection{}

Let us notice the the following: the set of all metamorphoses that is described by the formula (\ref{1.4}) is connected with the critical line
$\sigma=\frac 12$. In this paper we introduce a new complicated oscillatory $Q$-system, namely the quotient of two multiforms of the type (\ref{1.2}). Secondly, we obtain
an infinite set of metamorphoses of this system into an infinite set of infinite series of the type
\be \label{1.5}
\left|\sum_{n=1}^\infty \frac{\mu(n)}{n^{\sigma+it}}\right|,\ \sigma>1 ,
\ee
where $\mu(n)$ is  the M\" obius function. Consequently, we see that the old system (chrysalis), before metamorphosis, is defined on critical line
$\sigma=\frac 12$, and the new form (butterfly), after metamorphosis, is defined in the region $\sigma>1$. That is, mentioned sets are disconnected each to other
(their distance reads $1/2$).

\section{Theorem}

\subsection{}

Now we introduce new complicated oscillatory system as a quotient of two multiforms of type (\ref{1.2}) (the Q-system)
\be \label{2.1}
\begin{split}
& G(x_1,\dots,x_k;y_1,\dots,y_k)=\prod_{r=1}^k\left|\frac{Z(x_r)}{Z(y_r)}\right|= \\
& = \prod_{r=1}^k
\left|
\frac
{\sum_{n\leq \tau(x_r)}\frac{2}{\sqrt{n}}\cos\{\vth(x_r)-x_r\ln n\}+R(x_r)}
{\sum_{n\leq \tau(y_r)}\frac{2}{\sqrt{n}}\cos\{\vth(y_r)-y_r\ln n\}+R(y_r)}
\right|, \\
& (x_1,\dots,x_k)\in M_k^1,\quad (y_1,\dots,y_k)\in M_k^2, \\
& R(t)=\mcal{O}(t^{-1/4}), \\
& k\leq k_0\in\mbb{N}
\end{split}
\ee
($k_0$ is an arbitrary and fixed number), where
\be \label{2.2}
\begin{split}
& M_k^1=\{(x_1,\dots,x_k)\in (T_0,+\infty)^k:\ T_0<x_1<x_2<\dots<x_k,  \\
&  x_r\not=\gamma:\  \zeta\left(\frac 12+i\gamma\right)=0,\ r=1,\dots,k  \}, \\
& M_k^2=\{(y_1,\dots,y_k)\in (T_0,+\infty)^k:\ T_0<y_1<y_2<\dots<y_k,  \\
&  y_r\not=\gamma:\  \zeta\left(\frac 12+i\gamma\right)=0,\ r=1,\dots,k  \} .
\end{split}
\ee

\subsection{}

The following theorem holds true about the set of metamorphoses of the Q-system (\ref{2.1})

\begin{mydef1}
Let
\be \label{2.3}
[T,T+U]\rightarrow [\overset{1}{T},\overset{1}{\wideparen{T+U}}] ,\dots, [\overset{k}{T},\overset{k}{\wideparen{T+U}}]
\ee
where
\bdis
\begin{split}
& [\overset{r}{T},\overset{r}{\wideparen{T+U}}],\ r=1,\dots,j,\ k\leq k_0 \\
& U=U(T,\Theta)=\ln\ln T+\Theta \ln T,\ \Theta\in [0,1]
\end{split}
\edis
be the $r$-th reversely iterated segment corresponding to the first segment in (\ref{2.3}). Let furthermore
\be \label{2.4}
\sigma\in [1+\epsilon,+\infty),\ \epsilon>0
\ee
where $\epsilon$ is sufficiently small and fixed number. Then there is a sufficiently big
\bdis
T_0=T_0(\epsilon)
\edis
such that for every $T>T_0$ and every admissible $\sigma,\Theta,k,\epsilon$ there are functions
\be \label{2.5}
\begin{split}
& \alpha_r=\alpha_r(\sigma,T,\Theta,k,\epsilon),\ r=0,1,\dots,k, \\
& \alpha_r\not=\gamma:\ \zeta\left(\frac 12+i\gamma\right)=0, \\
& \beta_r=\beta_r(T,\Theta,k),\ r=1,\dots,k, \\
& \beta_r\not=\gamma:\ \zeta\left(\frac 12+i\gamma\right)=0,
\end{split}
\ee
such that
\be \label{2.6}
\begin{split}
& \prod_{r=1}^k
\left|
\frac
{\sum_{n\leq \tau(\alpha_r)}\frac{2}{\sqrt{n}}\cos\{\vth(\alpha_r)-\alpha_r\ln n\}+R(\alpha_r)}
{\sum_{n\leq \tau(\beta_r)}\frac{2}{\sqrt{n}}\cos\{\vth(\beta_r)-\beta_r\ln n\}+R(\beta_r)}
\right| \sim \\
& \sim \sqrt{\zeta(2\sigma)}\left|\sum_{n=1}^\infty \frac{\mu(n)}{n^{\sigma+i\alpha_0}}\right|,\ T\to\infty,
\end{split}
\ee
where $\mu(n)$ is the M\" obius function. Moreover, the sequences
\bdis
\{\alpha_r\}_{r=0}^k,\ \{\beta_r\}_{r=1}^k
\edis
have the following properties
\be \label{2.7}
\left\{
\begin{split}
& T<\alpha_0<\alpha_1<\dots<\alpha_k , \\
& T<\beta_1<\beta_2<\dots <\beta_k, \\
& \alpha_0\in (T,T+U),\ \alpha_r,\beta_r\in (\overset{r}{T},\overset{r}{\wideparen{T+U}}),\ r=1,\dots,k,
\end{split}
\right.
\ee
\be \label{2.8}
\left\{
\begin{split}
& \alpha_{r+1}-\alpha_r\sim (1-c)\pi(T),\ r=0,1,\dots,k-1 , \\
& \beta_{r+1}-\beta_r\sim (1-c)\pi(T),\ r=1,\dots,k-1 ,
\end{split}
\right.
\ee
where
\bdis
\pi(T)\sim \frac{T}{\ln T},\ T\to\infty
\edis
is the prime-counting function and $c$ is the Euler's constant.
\end{mydef1}

\subsection{}

Now, we give the following remarks.

\begin{remark}
The asymptotic behavior of the following sets
\be \label{2.9}
\{\alpha_r\}_{r=0}^k,\ \{\beta_r\}_{r=1}^k
\ee
is as follows (see (\ref{2.8}): if $T\to\infty$ then the points of every set in (\ref{2.9}) recede unboundedly each from other and all these points together
recede to infinity. Hence, at $T\to\infty$ each set of (\ref{2.9}) behaves as one-dimensional Friedmann-Hubble universe.
\end{remark}

\begin{remark}
In this Theorem we have obtained three resp. two parametric sets of control functions (=Golem's shem) for fixed
and admissible $k,\epsilon$
\be \label{2.10}
\begin{split}
 & \{ \alpha_0(\sigma,T,\Theta),\alpha_1(\sigma,T,\Theta),\dots,\alpha_k(\sigma,T,\Theta)\}, \\
 & \{ \beta_1(T,\Theta),\dots,\beta_k(T,\Theta)\}, \\
 & \sigma\in [1+\epsilon,+\infty),\ T\in (T_0,+\infty),\ \Theta\in [0,1],
\end{split}
\ee
of the metamorphoses (\ref{2.6}), (comp. Remark 1 and \cite{6}).
\end{remark}

\begin{remark}
The mechanism of metamorphosis is as follows. Let (comp. (\ref{2.2}), (\ref{2.10}))
\be \label{2.11}
\begin{split}
 & M_k^3=\{ \alpha_1(\sigma,T,\Theta),\dots,\alpha_k(\sigma,T,\Theta)\}, \\
 & M_k^4=\{ \beta_1(T,\Theta),\dots,\beta_k(T,\Theta)\}, \\
 & \sigma\in [1+\epsilon,+\infty),\ T\in (T_0,+\infty),\ \Theta\in [0,1],
\end{split}
\ee
where, of course,
\be \label{2.12}
\begin{split}
 & M_k^3\subset M_k^1\subset (T_0,+\infty)^k, \\
 & M_k^4\subset M_k^2\subset (T_0,+\infty)^k.
\end{split}
\ee
Now, if we obtain after random sampling (say) of the points
\bdis
(x_1,\dots,x_k),\ (y_1,\dots,y_k)
\edis
(see conditions (\ref{2.2}) on these) that
\be \label{2.13}
\begin{split}
 & (x_1,\dots,x_k)=(\alpha_1(\sigma,T,\Theta),\dots,\alpha_k(\sigma,T,\Theta))\in M_k^3(\epsilon), \\
 & (y_1,\dots,y_k)=(\beta_1(T,\Theta),\dots,\beta_k(T,\Theta))\in M_k^4
\end{split}
\ee
(see (\ref{2.11}), (\ref{2.12})) then - at the points (\ref{2.13}) - change (see (\ref{2.6})) the Q-system (\ref{2.1})
its old form (=chrysalis) into a new form (=butterfly) and the last is controlled by the function
$\alpha_0(\sigma,T,\Theta)$.
\end{remark}

\section{Riemann's oscillators as a basis of the Q-system (\ref{2.1})}

\subsection{}

The following local variant of the Riemann-Siegel formula holds true.

\begin{mydefSF}
\be \label{3.1}
\begin{split}
 & Z(t)=2\sum_{n\leq \tau(x_r)}\frac{1}{\sqrt{n}}\cos\left\{ t\ln\frac{\tau(x_r)}{n}-\frac{\tau(x_r)}{2}-\frac{\pi}{8}\right\}+ \\
 & +\mcal{O}(x_r^{-1/4}),\quad \tau(x_r)=\sqrt{\frac{x_r}{2\pi}}, \\
 & t\in [x_r,x_r+H],\ H\in (0,\sqrt[4]{x_r}],
\end{split}
\ee
where (comp. (\ref{2.2}))
\bdis
T_0<x_r,\ r=1,\dots,k.
\edis
\end{mydefSF}

\begin{proof}
 We will make the following transformations of the Riemann-Siegel formula
\be \label{3.2}
Z(t)=2\sum_{n\leq\tau(t)}\frac{1}{\sqrt{n}}\cos\{\vth(t)-t\ln n\}+\mcal{O}(t^{-1/4}).
\ee
\begin{itemize}
 \item[(a)] Since
 \bdis
 \begin{split}
  & \sum_{\tau(x_r)<n\leq \tau(x_r+H)}1=\mcal{O}\{\tau(x_r+H)-\tau(x_r)\}= \\
  & = \mcal{O}\left(\sqrt{\frac{x_r+H}{2\pi}}-\sqrt{\frac{x_r}{2\pi}}\right)=\mcal{O}\left(\frac{H}{\sqrt{x_r}}\right),
 \end{split}
 \edis
 then
 \bdis
 \begin{split}
  & \sum_{\tau(x_r)<n\leq \tau(x_r+H)}\frac{2}{\sqrt{n}}\cos\{\vth(t)-t\ln n\}= \\
  & = \mcal{O}\left(\frac{H}{\sqrt{x_r}}\frac{1}{\sqrt{\tau(x_r)}}\right)=
  \mcal{O}\left(\frac{H}{x_r^{3/2}}\right)=\mcal{O}(x_r^{-1/2}).
 \end{split}
 \edis
 Consequently, we have (see (\ref{3.2})) that
 \be \label{3.3}
 \begin{split}
  & Z(t)=2\sum_{n\leq \tau(x_r)}\frac{1}{\sqrt{n}}\cos\{\vth(t)-t\ln n\}+\mcal{O}(x_r^{-1/4}), \\
  & t\in [x_r,x_r+H].
 \end{split}
 \ee
\item[(b)] Next, we use the following formula
\be \label{3.4}
\begin{split}
 & \vth(t)=\vth(x_r)+\vth'(x_r)(t-x_r)+\frac 12\vth''(\xi_r)(t-x_r)^2, \\
 & \xi_r\in (x_r,x_r+H),\ t\in [x_r,x_r+H].
\end{split}
\ee
Since (see \cite{8}, pp. 221, 329)
\bdis
\begin{split}
 & \vth(t)=\frac t2\ln\frac{t}{2\pi}-\frac t2-\frac{\pi}{8}+\mcal{O}(t^{-1}), \\
 & \vth'(t)=\frac 12\ln\frac{t}{2\pi}+\mcal{O}(t^{-1}), \\
 & \vth''(t)\sim \frac{1}{2t},
\end{split}
\edis
then
\be \label{3.5}
\vth(x_r)=x_r\ln\tau(x_r)-\frac{x_r}{2}-\frac{\pi}{8}+\mcal{O}\left(\frac{1}{x_r}\right),
\ee
and
\be \label{3.6}
\begin{split}
 & \vth'(x_r)(t-x_r)=\left\{\frac 12\ln\frac{x_r}{2\pi}+\mcal{O}\left(\frac{1}{x_r}\right)\right\}(t-x_r)= \\
 & = \left\{\ln\tau(x_r)+\mcal{O}\left(\frac{1}{x_r}\right)\right\}(t-x_r)= \\
 & = t\ln\tau(x_r)-x_r\ln\tau(x_r)+\mcal{O}\left(\frac{t-x_r}{x_r}\right)= \\
 & = t\ln\tau(x_r)-x_r\ln\tau(x_r)+\mcal{O}\left(\frac{H}{x_r}\right).
\end{split}
\ee
Hence (see (\ref{3.3})--(\ref{3.6}))
\be\label{3.7}
\begin{split}
 & \vth(t)-t\ln n=t\ln\frac{\tau(x_r)}{n}-\frac{x_r}{2}-\frac{\pi}{8}+\mcal{O}\left(\frac{1+H}{x_r}\right)= \\
 & = t\ln\frac{\tau(x_r)}{n}-\frac{x_r}{2}-\frac{\pi}{8}+\mcal{O}(x_r^{-1/4}).
\end{split}
\ee
Putting (\ref{3.7}) into (\ref{3.3}) we obtain the result (\ref{3.1}).
\end{itemize}
\end{proof}

\subsection{}

It is natural to introduce the following terminology (based on the formulae (\ref{3.1}),(\ref{3.2}))

\begin{remark}
We shall call:
\begin{itemize}
\item[(a)] the formula (\ref{3.1}) as the local spectral form of the Riemann-Siegel formula (\ref{3.2})
\item[(b)] the sets
\bdis
\begin{split}
 & \{\omega_n(x_r)\}_{n\leq\tau(x_r)},\ \omega_n(x_r)=\ln\frac{\tau(x_r)}{n}, \\
 & \{\omega_n(y_r)\}_{n\leq\tau(y_r)},\ \omega_n(y_r)=\ln\frac{\tau(y_r)}{n}, \\
 & r=1,\dots,k
\end{split}
\edis
as the local spectrum of the cyclic frequencies of the Q-system (\ref{2.1}),
\item[(c)] the expressions
\be \label{3.8}
\begin{split}
 & \frac{2}{\sqrt{n}}\cos\left\{ t\omega_n(x_r)-\frac{x_r}{2}-\frac{\pi}{8}\right\},\dots \\
 & t\in [x_r,x_r+H],\ 1\leq n\leq \tau(x_r),\dots
\end{split}
\ee
as the local Riemann's oscillators with the set of incoherent local phase \emph{constants}
\bdis
\left\{ -\frac{x_r}{2}-\frac{\pi}{8}\right\}, \left\{ -\frac{y_r}{2}-\frac{\pi}{8}\right\},
\edis
and, moreover, with the set of non-synchronized local times $t=t(x_r)$.
\end{itemize}
\end{remark}

Now, based on our Remark 5 we give the following

\begin{remark}
We see that the Q-system
\bdis
G(x_1,\dots,x_k;y_1,\dots,y_k)
\edis
(see (\ref{2.1})) expresses the complicated oscillating process that is generated by oscillations of a big number
of the local Riemann's oscillators (\ref{3.8}). Just for this oscillating Q-system we have the oscillators (\ref{3.8}).
Just for this oscillating Q-system we have obtained the infinite set of metamorphoses described by the formula (\ref{2.6}).
For example, if
\bdis
k=F_4=2^{2^4}+1=65537
\edis
where $F_4$ is the Fermat-Gauss prime, then the above mentioned big number of Riemann's oscillators (interacting oscillators)
is bigger than
\bdis
2\left(\sqrt{\frac{T_0}{2\pi}}\right)^{65537}.
\edis
\end{remark}

\section{Proof of the Theorem}

\subsection{}

Let us remind the formula (see \cite{4}, (2.1))
\be \label{4.1}
\int_T^{T+U} |\zeta(\sigma+it)|^2{\rm d}t=\zeta(2\sigma)U+\mcal{O}(1).
\ee
This formula holds true uniformly for
\bdis
T,U>0,\ \sigma\geq 1+\epsilon,\ \epsilon>0,
\edis
where $\epsilon$ is arbitrary small fixed number and the $\mcal{O}$-constant depends on the choice of that $\epsilon$. Since
\bdis
\zeta(2\sigma)U+\mcal{O}(1)=\zeta(2\sigma)U\left\{1+\mcal{O}\left(\frac 1U\right)\right\}
\edis
then the formula (\ref{4.1}) is asymptotic one. For example, in the case
\bdis
U\geq \ln\ln T,\ \sigma\in[1+\epsilon,+\infty).
\edis
We use in this paper the following local version of the formula (\ref{4.1})
\be \label{4.2}
\begin{split}
& \int_T^{T+U(T,\Theta)}|\zeta(\sigma+it)|^2{\rm d}t\sim \zeta(2\sigma)U(T,\Theta), \\
& U(T,\Theta)=\ln\ln T+\Theta\ln T, \\
& \sigma\in[1+\epsilon,+\infty),\ \Theta\in [0,1].
\end{split}
\ee

\subsection{}

Since (see (\ref{4.2}))
\be \label{4.3}
U(T,\Theta)=o\left(\frac{T}{\ln T}\right),
\ee
then we have by our lemma (see \cite{5}, (7.1), (7.2)) that
\be \label{4.4}
\begin{split}
& \int_T^{T+U(T,\Theta)}|\zeta(\sigma+it)|^2{\rm d}t= \\
& =\int_{\overset{k}{T}}^{\overset{k}{\wideparen{T+U}}}|\zeta[\sigma+i\vp_1^k(t)]|^2\prod_{r=0}^{k-1}\tilde{Z}^2[\vp_1^r(t)]{\rm d}t
\end{split}
\ee
i.e. (see (\ref{4.2}), (\ref{4.4}))
\be \label{4.5}
\int_{\overset{k}{T}}^{\overset{k}{\wideparen{T+U}}}|\zeta[\sigma+i\vp_1^k(t)]|^2\prod_{r=0}^{k-1}\tilde{Z}^2[\vp_1^r(t)]{\rm d}t
\sim \zeta(2\sigma)U.
\ee
Next, we obtain by using the mean-value theorem in (\ref{4.5}) that
\be \label{4.6}
\begin{split}
& \prod_{r=0}^{k-1} \tilde{Z}^2[\vp_1^r(t)]\sim \zeta(2\sigma)
\frac{U}{\overset{k}{\wideparen{T+U}}-\overset{k}{T}}
\frac{1}{|\zeta[\sigma+i\vp_1^k(d)]|^2}, \\
& d=d(\sigma,T,\Theta,k,\epsilon)\in (\overset{k}{T},\overset{k}{\wideparen{T+U}}) .
\end{split}
\ee

\subsection{}

Further, we have (comp. \cite{5}, Property 2, (6.4) and \cite{6}, (5.4)) that
\bdis
d\in (\overset{k}{T},\overset{k}{\wideparen{T+U}}) \ \Rightarrow \
\vp_1^r(d)\in (\overset{k-r}{T},\overset{k-r}{\wideparen{T+U}}),\ r=0,1,\dots,k,
\edis
i.e.
\be \label{4.7}
\begin{split}
& \vp_1^0(d)=\alpha_k\in (\overset{k}{T},\overset{k}{\wideparen{T+U}}) , \\
& \vp_1^1(d)=\alpha_{k-1}\in (\overset{k-1}{T},\overset{k-1}{\wideparen{T+U}}) , \\
& \vdots \\
& \vp_1^{k-2}(d)=\alpha_{2}\in (\overset{2}{T},\overset{2}{\wideparen{T+U}}) , \\
& \vp_1^{k-1}(d)=\alpha_{1}\in (\overset{1}{T},\overset{1}{\wideparen{T+U}}) , \\
& \vp_1^k(d)=\alpha_0\in (T,T+U).
\end{split}
\ee
Consequently, from (\ref{4.6}) by (\ref{4.7}) the formula
\be \label{4.8}
\begin{split}
& \prod_{l=1}^k \tilde{Z}^2(\alpha_l)\sim \zeta(2\sigma)\frac{U}{\overset{k}{\wideparen{T+U}}-\overset{k}{T}}
\frac{1}{|\zeta(\sigma+i\alpha_0)|^2}, \\
& \alpha_r=\alpha_r(\sigma,T,\Theta,k,\epsilon),\ r=0,1,\dots,k
\end{split}
\ee
follows.

\begin{remark}
Wee see that corresponding inclusions for $\alpha_r$ in (\ref{2.7}) follows from (\ref{4.7}).
\end{remark}

\subsection{}

Next, we introduce the infinite collection of the following disconnected sets
\be \label{4.9}
\Delta(T,\Theta,k)=\bigcup_{r=0}^k [\overset{r}{T},\overset{r}{\wideparen{T+U(T,\Theta)}}],
\ee
(comp. (\ref{4.2})). Let us remind the following properties of (\ref{4.9}) (see (\ref{4.3}), comp. \cite{5}, (2.5) -- (2.7), (2.9)): since
\bdis
U(t,\Theta)=o\left(\frac{T}{\ln T}\right)
\edis
then
\be \label{4.10}
\begin{split}
& |[\overset{r}{T},\overset{r}{\wideparen{T+U}}]|=\overset{r}{\wideparen{T+U}}-\overset{r}{T}=o\left(\frac{T}{\ln T}\right),\ r=1,\dots,k, \\
& |[\overset{r-1}{\wideparen{T+U}},\overset{r}{T}]|\sim (1-c)\pi(T) , \\
& [T,T+U]\prec [\overset{1}{T},\overset{1}{\wideparen{T+U}}]\prec \dots\prec [\overset{k}{T},\overset{k}{\wideparen{T+U}}].
\end{split}
\ee
Consequently, the properties (\ref{2.7}), (\ref{2.8}) follows from (\ref{4.7}), (\ref{4.10}) immediately.

\subsection{}

Next, we have (see \cite{2}, (3.9), \cite{3}, (9.1), (9.2), comp. \cite{6}, (4.1), (4.2)) that
\bdis
\tilde{Z}^2(t)=\frac{{\rm d}\vp_1(t)}{{\rm d}t},\ \vp_1(t)=\frac 12\vp(t),
\edis
where
\be \label{4.11}
\begin{split}
& \tilde{Z}^2(t)=\frac{Z^2(t)}{2\Phi'_\vp[\vp(t)]}=\frac{\left|\zf\right|^2}{\omega(t)}, \\
& \omega(t)=\left\{ 1+\mcal{O}\left(\frac{\ln\ln t}{\ln t}\right)\right\}\ln t.
\end{split}
\ee
We call the function $\vp_1(t)$ the Jacob's ladder (see our papers \cite{2}, \cite{3}) according to Jacob's
dream in Chumash, Bereishis, 28:12. Further we have (comp. \cite{5}, (4.3))
\be \label{4.12}
\ln t\sim \ln T,\ \forall t\in (T,\overset{k}{\wideparen{T+U}}),
\ee
i.e. we have (see (\ref{4.8}), (\ref{4.11}), (\ref{4.12})) that
\be \label{4.13}
\tilde{Z}^2(\alpha_r)\sim \frac{Z^2(\alpha_r)}{\ln T},\ T\to\infty.
\ee
Consequently, the following formula
\be \label{4.14}
\begin{split}
 & \prod_{r=1}^k Z^2(\alpha_r)\sim \\
 & \sim \zeta(2\sigma)\frac{U}{\overset{k}{\wideparen{T+U}}-\overset{k}{T}}\ln^kT
 \frac{1}{|\zeta(\sigma+i\alpha_0)|^2}
\end{split}
\ee
(see (\ref{4.8}), (\ref{4.13})) holds true.

\subsection{}

Next, let us remind the following formula (see conditions (\ref{4.3}) and \cite{3}, (9.5), comp.
\cite{5}, (7.1), (7.2))
\bdis
\int_T^{T+U}f(t){\rm d}t=
\int_{\overset{k}{T}}^{\overset{k}{\wideparen{T+U}}}f[\vp_1^k(t)]
\prod_{r=0}^{k-1} \tilde{Z}^2[\vp_1^r(t)]{\rm d}t
\edis
holds true. If we put
\bdis
f(t)=1
\edis
then we obtain
\be \label{4.15}
U=\int_{\overset{k}{T}}^{\overset{k}{\wideparen{T+U}}}\prod_{r=0}^{k-1} \tilde{Z}^2[\vp_1^r(t)]{\rm d}t.
\ee
Further, we have by making use of the mean-value theorem in (\ref{4.15})
\be \label{4.16}
\begin{split}
 & \prod_{r=0}^{k-1} \tilde{Z}^2[\vp_1^r(e)]=\frac{U}{\overset{k}{\wideparen{T+U}}-\overset{k}{T}}, \\
 & e=e(T,\Theta,k)\in (\overset{k}{T},\overset{k}{\wideparen{T+U}}).
\end{split}
\ee
Since
\bdis
e\in (\overset{k}{T},\overset{k}{\wideparen{T+U}}) \ \Rightarrow \
\vp_1^r(e)\in (\overset{k-r}{T},\overset{k-r}{\wideparen{T+U}}),\ r=0,1,\dots,k-1,
\edis
then we have (similarly to (\ref{4.7})) that
\be \label{4.17}
\begin{split}
 & \vp_1^0(e)=\beta_k\in (\overset{k}{T},\overset{k}{\wideparen{T+U}}), \\
 & \vp_1^1(e)=\beta_{k-1}\in (\overset{k-1}{T},\overset{k-1}{\wideparen{T+U}}), \\
 & \vdots \\
 & \vp_1^{k-2}(e)=\beta_2\in (\overset{2}{T},\overset{2}{\wideparen{T+U}}), \\
 & \vp_1^{k-1}(e)=\beta_1\in (\overset{1}{T},\overset{1}{\wideparen{T+U}}) .
\end{split}
\ee
Now, we obtain from (\ref{4.16}) by (\ref{4.17}), (comp. (\ref{4.13})) the following formula
\be \label{4.18}
\begin{split}
 & \prod_{r=1}^k Z^2(\beta_r)\sim \frac{U}{\overset{k}{\wideparen{T+U}}-\overset{k}{T}}\ln^kT,\ T\to\infty, \\
 & \beta_r=\beta_r(T,\Theta,k),\ r=1,\dots,k.
\end{split}
\ee

\begin{remark}
Of course, our sequence $\{\beta_r\}$, similarly to the sequence $\{\alpha_r\}$ has the properties
listen in (\ref{2.7}), (\ref{2.8}).
\end{remark}

Consequently, the following factorization formula (comp. \cite{6}), (see (\ref{4.14}), (\ref{4.18}))
\be \label{4.19}
\prod_{r=1}^k \left|\frac{Z[\alpha_r(\sigma)]}{Z(\beta_r)}\right|\sim
\frac{\sqrt{\zeta(2\sigma)}}{|\zeta[\sigma+i\alpha_0(\sigma)]|},\ T\to\infty
\ee
holds true. Finally, the formula (\ref{2.6}) follows from (\ref{4.19}) (see (\ref{1.5}), (\ref{2.1})).

\thanks{I would like to thank Michal Demetrian for his help with electronic version of this paper.}

\end{document}